\documentclass[11pt,a4paper]{amsart}

\usepackage{a4wide}
\usepackage[T1]{fontenc}
\usepackage[british]{babel}
\usepackage{lmodern}
\usepackage{microtype}
\usepackage{mathtools,amssymb,amsthm}
\usepackage{enumitem}
\usepackage{url}
\usepackage[hidelinks]{hyperref}
\hypersetup{
  pdftitle={Higher traces as boundary averages on finite-dimensional normed spaces},
  pdfauthor={Tomasz Kania},
  pdfkeywords={higher traces, exterior powers, cone measure, norming functionals, hypersurface measure, Auerbach basis, spherical harmonics, Minkowski identity}
}
\usepackage[nameinlink,capitalise,noabbrev]{cleveref}

\numberwithin{equation}{section}
\theoremstyle{plain}
\newtheorem{theorem}{Theorem}[section]
\newtheorem{proposition}[theorem]{Proposition}
\newtheorem{lemma}[theorem]{Lemma}
\newtheorem{corollary}[theorem]{Corollary}

\newtheorem{maintheorem}{Theorem}

\theoremstyle{definition}

\newtheorem{example}[theorem]{Example}

\theoremstyle{remark}

\DeclareMathOperator{\End}{End}

\DeclareMathOperator{\rank}{rank}

\DeclareMathOperator{\Sym}{Sym}
\DeclareMathOperator{\tr}{tr}
\DeclareMathOperator{\vol}{vol}
\let\div\relax
\DeclareMathOperator{\div}{div}

\newcommand{\R}{\mathbb{R}}
\newcommand{\dd}{\,\mathrm{d}}
\newcommand{\Id}{I}
\newcommand{\sphere}{\mathbb{S}}
\newcommand{\ip}[2]{\left\langle #1,#2\right\rangle_2}
\newcommand{\abs}[1]{\left|#1\right|}
\newcommand{\norm}[1]{\left\|#1\right\|}
\newcommand{\restr}{\mathbin{|}}

\title[Higher traces as boundary averages]{Higher traces as boundary averages on finite-dimensional normed spaces}
\author{Tomasz Kania}
\address[T.~Kania]{Mathematical Institute\\Czech Academy of Sciences\\\v Zitn\'a 25 \\115 67 Praha 1\\Czech Republic  and  Institute of Mathematics and Computer Science\\ Jagiellonian University\\ {\L}ojasiewicza 6, 30-348 Krak\'{o}w, Poland
}
\email{kania@math.cas.cz, tomasz.marcin.kania@gmail.com}
\thanks{RVO: 67985840.}
\date{}


\begin{document}

\begin{abstract}
Let $X$ be an $N$-dimensional real normed space, let $1\leqslant k\leqslant N$, and set $V=\Lambda^kX$ and $m=\binom Nk$. We characterise the probability measures $\eta$ on the unit sphere of $V$ for which
\[
\tr(\Lambda^kA)=m\int w^\sharp\big((\Lambda^kA)w\big)\,\dd\eta(w)
\]
holds for every $A\in\End(X)$: this is equivalent to $m\int w\otimes w^\sharp\,\dd\eta(w)=\Id_V$. The cone probability measure always satisfies this condition, giving a canonical higher-trace formula for every norm. Normalised Euclidean hypersurface measure also does so under a scalar-commutant symmetry hypothesis, including spaces with a $1$-symmetric basis. We further obtain atomic and polyhedral formulae and show that, within a natural power-weighted family, cone measure is the unique universally isotropic member; for hypersurface measure the first-order obstruction is precisely the degree-$2$ spherical harmonic component of the support function.
\end{abstract}

\maketitle

\section{Introduction}

Let $X$ be an $N$-dimensional real vector space and let $A\in\End(X)$. For $0\leqslant k\leqslant N$, the $k$th higher trace of $A$ is
\[
\lambda_k(A):=\tr(\Lambda^kA),
\]
where $\lambda_0(A)=1$. Equivalently,
\begin{equation}\label{eq:char-poly}
\det(\Id_X-tA)=\sum_{k=0}^N(-1)^k\lambda_k(A)t^k.
\end{equation}
Thus $\lambda_1(A)=\tr(A)$, $\lambda_N(A)=\det(A)$, and $\lambda_k(A)$ is the sum of the principal $k\times k$ minors of a matrix representing $A$.

When $X=\R^N$ is Euclidean, Eberlein's formula expresses $\lambda_k(A)$ as an average over the Grassmannian:
\begin{equation}\label{eq:eberlein-intro}
\lambda_k(A)=\binom Nk\int_{G_{k,N}}\det\big(P_EA\restr_E\big)\,\dd\sigma_{k,N}(E),
\end{equation}
where $\sigma_{k,N}$ is the $O(N)$-invariant probability measure; see \cite{Eberlein,MorrisonHigher}. Formula \eqref{eq:eberlein-intro} is the exterior-power analogue of the familiar spherical formula for the trace. The latter extends to several classes of non-Euclidean norms, notably to spaces with a $1$-symmetric basis \cite{KaniaTrace,KM}.

The purpose of this paper is to separate the algebraic and geometric ingredients of such formulae. Let $X$ now be normed, let
\[
V:=\Lambda^kX,
\qquad
m:=\dim V=\binom Nk,
\]
and equip $V$ with the projective exterior norm. For $w\in S_V$, a norming functional is a functional $w^\sharp\in S_{V^*}$ such that $w^\sharp(w)=1$. Given a probability measure $\eta$ on $S_V$ and a measurable choice of norming functionals, define
\begin{equation}\label{eq:T-intro}
T_\eta:=m\int_{S_V}w\otimes w^\sharp\,\dd\eta(w),
\qquad
(w\otimes w^\sharp)(z):=w^\sharp(z)w.
\end{equation}
The operator $T_\eta$ is the appropriate isotropy tensor. Our first result gives an exact criterion, including the converse for higher traces.

\begin{maintheorem}\label{thm:intro-A}
Let $1\leqslant k\leqslant N$. For a probability measure $\eta$ and a measurable norming selection on $S_{\Lambda^kX}$, the following are equivalent:
\begin{enumerate}[label=\textup{(\roman*)}]
\item $T_\eta=\Id_{\Lambda^kX}$;
\item for every $B\in\End(\Lambda^kX)$,
\[
\tr(B)=m\int w^\sharp(Bw)\,\dd\eta(w);
\]
\item for every $A\in\End(X)$,
\[
\lambda_k(A)=m\int w^\sharp\big((\Lambda^kA)w\big)\,\dd\eta(w).
\]
\end{enumerate}
\end{maintheorem}

The implication from (iii) to (i) rests on a simple but useful fact: relative to the standard exterior basis, every matrix unit on $\Lambda^kX$ is equal to $\Lambda^kA$ for an explicitly chosen rank-$k$ operator $A$ on $X$. Thus the family $\{\Lambda^kA:A\in\End(X)\}$ linearly spans the whole algebra $\End(\Lambda^kX)$.

There is always a discrete isotropic measure: an Auerbach basis $(w_j,w_j^\sharp)_{j=1}^m$ of $V$ gives
\[
\Id_V=\sum_{j=1}^m w_j\otimes w_j^\sharp.
\]
The corresponding probability measure has exactly $m$ atoms, and this number is minimal. The main geometric point is that there is also a canonical continuous choice, namely cone measure.

Fix an auxiliary Euclidean structure on $V$, let $B_V$ be its norm ball, and write $n_V(w)$ for the Euclidean outer unit normal at a regular point $w\in\partial B_V$. The cone probability measure is
\[
\nu_V(E)=\frac{\vol\{tw:w\in E,\ 0\leqslant t\leqslant1\}}{\vol(B_V)}.
\]
It does not depend on the normalisation of Lebesgue measure and is natural under linear isomorphisms.

\begin{maintheorem}\label{thm:intro-B}
For every finite-dimensional normed space $V$, cone measure satisfies
\[
(\dim V)\int_{S_V}w\otimes w^\sharp\,\dd\nu_V(w)=\Id_V.
\]
Consequently, for $V=\Lambda^kX$ and every $A\in\End(X)$,
\begin{align*}
\lambda_k(A)
&=m\int_{S_V}w^\sharp\big((\Lambda^kA)w\big)\,\dd\nu_V(w)\\
&=\frac{1}{\vol(B_V)}
  \int_{\partial B_V}\ip{(\Lambda^kA)w}{n_V(w)}\,\dd\mathcal H^{m-1}(w).
\end{align*}
\end{maintheorem}

The proof is the Gauss--Green theorem applied to linear vector fields, equivalently the matrix-valued Minkowski identity
\[
\int_{\partial B_V}w\otimes n_V(w)\,\dd\mathcal H^{m-1}(w)
=\vol(B_V)\Id_V.
\]
No smoothness or symmetry is required.

The normalised Euclidean hypersurface measure $\mu_V$ is less canonical: it depends on the auxiliary Euclidean structure and it need not be isotropic. Symmetry nevertheless gives a clean sufficient condition.

\begin{maintheorem}\label{thm:intro-C}
Let $G$ be a compact group of norm isometries of $V$, and choose an auxiliary Euclidean inner product for which $G$ acts orthogonally. If
\[
\End(V)^G:=\{R\in\End(V):RQ=QR\text{ for every }Q\in G\}=\R\Id_V,
\]
then $T_{\mu_V}=\Id_V$. In particular, if $X$ has a $1$-symmetric basis, then the hypersurface formula in Theorem~\textup{A} holds on $\Lambda^kX$ for every $1\leqslant k\leqslant N$.
\end{maintheorem}

Finally, we compare hypersurface and cone measure within a natural one-parameter family. Let $K\subset\R^m$ be a smooth strictly convex body containing the origin, let $h_K$ be its support function, and define a probability measure on $\partial K$ by a density proportional to $h_K^\alpha$ with respect to surface area. The cases $\alpha=0$ and $\alpha=1$ are hypersurface and cone measure, respectively.

\begin{maintheorem}\label{thm:intro-D}
Let $m\geqslant2$, let $g\in C^\infty(\sphere^{m-1})$, and let $K_\varepsilon$ have support function $h_\varepsilon=1+\varepsilon g$ for sufficiently small $\varepsilon$. If $T_\alpha(\varepsilon)$ is the isotropy tensor associated with the probability measure whose density is proportional to $h_\varepsilon^\alpha$, then
\begin{equation}\label{eq:intro-first-variation}
\left.\frac{\dd}{\dd\varepsilon}\right|_{\varepsilon=0}T_\alpha(\varepsilon)
=m(\alpha-1)\int_{\sphere^{m-1}}g(u)
\left(u\otimes u-\frac1m\Id_m\right)\,\dd\omega(u),
\end{equation}
where $\omega$ is rotationally invariant probability measure. Hence only the degree-$2$ spherical harmonic component of $g$ contributes. If
\[
P_2g(u)=u^\top Su,
\qquad S\in\Sym(\R^m),\quad \tr(S)=0,
\]
then the derivative in \eqref{eq:intro-first-variation} equals
\[
\frac{2(\alpha-1)}{m+2}S.
\]
In particular, $\alpha=1$ is the unique exponent for which this family is isotropic for every smooth centrally symmetric convex body.
\end{maintheorem}

The paper is organised as follows. In \cref{sec:isotropy} we develop the general isotropy formalism and the minimal atomic formula. In \cref{sec:exterior} we prove the exact higher-trace criterion and recall the Grassmannian formula. Cone measure and polyhedral formulae are treated in \cref{sec:cone}. The symmetry criterion for hypersurface measure is proved in \cref{sec:symmetry}. The first-variation theorem and explicit counterexamples occupy \cref{sec:variation}.

\section{Isotropy operators and atomic formulae}\label{sec:isotropy}

Let $V$ be an $m$-dimensional real normed space. Write $B_V$ and $S_V$ for its closed unit ball and unit sphere. Suppose that $\eta$ is a probability measure on $S_V$ and that
\[
w\longmapsto w^\sharp\in S_{V^*},
\qquad
w^\sharp(w)=1,
\]
is an $\eta$-measurable choice of norming functionals. Define $T_\eta$ by
\begin{equation}\label{eq:T-general}
T_\eta=m\int_{S_V}w\otimes w^\sharp\,\dd\eta(w).
\end{equation}
The integral is unproblematic because the integrand is bounded and takes values in a finite-dimensional space.

\begin{lemma}\label{lem:trace-pairing}
For every $B\in\End(V)$,
\begin{equation}\label{eq:trace-pairing}
m\int_{S_V}w^\sharp(Bw)\,\dd\eta(w)=\tr(BT_\eta).
\end{equation}
Consequently,
\[
\tr(B)=m\int_{S_V}w^\sharp(Bw)\,\dd\eta(w)
\quad\text{for every }B\in\End(V)
\]
if and only if $T_\eta=\Id_V$.
\end{lemma}

\begin{proof}
For $w\in V$ and $f\in V^*$,
\[
\tr\big(B(w\otimes f)\big)=f(Bw).
\]
Integrating gives \eqref{eq:trace-pairing}. The trace pairing $(B,C)\mapsto\tr(BC)$ on $\End(V)$ is non-degenerate, so $\tr(B(T_\eta-\Id_V))=0$ for every $B$ if and only if $T_\eta=\Id_V$.
\end{proof}

The trace of every isotropy tensor is already correctly normalised:
\begin{equation}\label{eq:trace-T}
\tr(T_\eta)=m\int_{S_V}w^\sharp(w)\,\dd\eta(w)=m.
\end{equation}
It is therefore natural to call
\[
\mathcal A_\eta:=T_\eta-\Id_V
\]
the anisotropy tensor of the pair $(\eta,w^\sharp)$. With an auxiliary Euclidean inner product on $V$, \cref{lem:trace-pairing} gives the quantitative estimate
\begin{equation}\label{eq:HS-bound}
\left|m\int w^\sharp(Bw)\,\dd\eta(w)-\tr(B)\right|
=\abs{\tr(B\mathcal A_\eta)}
\leqslant\norm{B}_{\mathrm{HS}}\norm{\mathcal A_\eta}_{\mathrm{HS}}.
\end{equation}

There is always an isotropic measure with minimal finite support.

\begin{lemma}\label{lem:auerbach}
Every $m$-dimensional normed space $V$ admits vectors $v_1,\ldots,v_m\in S_V$ and functionals $v_1^\sharp,\ldots,v_m^\sharp\in S_{V^*}$ such that
\[
v_i^\sharp(v_j)=\delta_{ij}.
\]
\end{lemma}

\begin{proof}
Fix a non-zero alternating $m$-linear form $D$ on $V$ and maximise
\[
\abs{D(v_1,\ldots,v_m)}
\]
over the compact set $B_V^m$. A maximising $m$-tuple is a basis, and each $v_i$ has norm one, since otherwise rescaling it would increase the determinant. Let $v_i^\sharp$ be the corresponding coordinate functionals. Replacing $v_i$ by any $v\in B_V$ shows that
\[
\abs{v_i^\sharp(v)}\leqslant1.
\]
Thus $\norm{v_i^\sharp}\leqslant1$, while $v_i^\sharp(v_i)=1$ gives equality.
\end{proof}

\begin{proposition}\label{prop:atomic}
Let $(v_i,v_i^\sharp)_{i=1}^m$ be as in \cref{lem:auerbach}. Then
\begin{equation}\label{eq:auerbach-identity}
\Id_V=\sum_{i=1}^m v_i\otimes v_i^\sharp,
\end{equation}
and the probability measure
\[
\eta_{\mathrm{A}}:=\frac1m\sum_{i=1}^m\delta_{v_i}
\]
is isotropic. Moreover, every isotropic atomic probability measure on $S_V$ has at least $m$ atoms. If it has exactly $m$ atoms, then all weights are $1/m$ and the atoms together with the selected norming functionals form an Auerbach basis.
\end{proposition}

\begin{proof}
Identity \eqref{eq:auerbach-identity} is the usual resolution of the identity in a basis and its biorthogonal functionals. Hence $T_{\eta_{\mathrm A}}=\Id_V$.

Now suppose
\[
\eta=\sum_{i=1}^r a_i\delta_{w_i},
\qquad a_i>0,
\qquad \sum_{i=1}^r a_i=1,
\]
is isotropic. Then
\[
\Id_V=m\sum_{i=1}^r a_iw_i\otimes w_i^\sharp,
\]
so $m=\rank(\Id_V)\leqslant r$. If $r=m$, the vectors $w_i$ form a basis. Let $g_i$ be their biorthogonal functionals. Applying $g_j$ to the displayed identity gives
\[
g_j=ma_jw_j^\sharp.
\]
Evaluation at $w_j$ yields $a_j=1/m$, and then $g_j=w_j^\sharp$.
\end{proof}

\section{Exterior powers and the exact higher-trace criterion}\label{sec:exterior}

Let $X$ be an $N$-dimensional real normed space and let $1\leqslant k\leqslant N$. On $V=\Lambda^kX$ we use the projective exterior norm
\begin{equation}\label{eq:exterior-projective}
\norm{w}_{\wedge,\pi}
=\inf\left\{\sum_{r=1}^s\prod_{j=1}^k\norm{x_j^{(r)}}:
 w=\sum_{r=1}^s x_1^{(r)}\wedge\cdots\wedge x_k^{(r)}\right\}.
\end{equation}
This is the quotient of the projective tensor norm under the canonical alternating map; see, for example, \cite{QuasThieullenZarrabi}.

\begin{lemma}\label{lem:exterior-bound}
For every $A\in\End(X)$,
\[
\norm{\Lambda^kA}\leqslant\norm{A}^k.
\]
In particular, if $A$ is a linear isometry of $X$, then $\Lambda^kA$ is a linear isometry of $(\Lambda^kX,\norm{\cdot}_{\wedge,\pi})$.
\end{lemma}

\begin{proof}
For every representation occurring in \eqref{eq:exterior-projective},
\begin{align*}
\norm{(\Lambda^kA)w}_{\wedge,\pi}
&\leqslant\sum_{r=1}^s\prod_{j=1}^k\norm{Ax_j^{(r)}}\\
&\leqslant\norm{A}^k\sum_{r=1}^s\prod_{j=1}^k\norm{x_j^{(r)}}.
\end{align*}
Taking the infimum proves the estimate. If $A$ is an isometry, apply the estimate to $A$ and $A^{-1}$.
\end{proof}

The next observation is responsible for the converse in Theorem~\ref{thm:intro-A}.

\begin{lemma}\label{lem:matrix-units}
The linear span of
\[
\{\Lambda^kA:A\in\End(X)\}
\]
is $\End(\Lambda^kX)$. More precisely, relative to any exterior basis, every matrix unit on $\Lambda^kX$ equals $\Lambda^kA$ for a suitable rank-$k$ operator $A$ on $X$.
\end{lemma}

\begin{proof}
Fix a basis $e_1,\ldots,e_N$ of $X$. For a $k$-element subset
\[
I=\{i_1<\cdots<i_k\},
\]
write $e_I=e_{i_1}\wedge\cdots\wedge e_{i_k}$. Given two $k$-element subsets $I$ and $J=\{j_1<\cdots<j_k\}$, define $A_{I,J}\in\End(X)$ by
\[
A_{I,J}e_{j_r}=e_{i_r}\quad(1\leqslant r\leqslant k),
\qquad
A_{I,J}e_\ell=0\quad(\ell\notin J).
\]
Then
\[
(\Lambda^kA_{I,J})e_J=e_I,
\qquad
(\Lambda^kA_{I,J})e_K=0\quad(K\neq J).
\]
Thus $\Lambda^kA_{I,J}$ is the matrix unit $e_I\otimes e_J^*$.
\end{proof}

\begin{theorem}\label{thm:exact-criterion}
Let $V=\Lambda^kX$, $m=\binom Nk$, and let $\eta$ be a probability measure on $S_V$ with a measurable norming selection. The following are equivalent:
\begin{enumerate}[label=\textup{(\roman*)}]
\item $T_\eta=\Id_V$;
\item $\tr(B)=m\int w^\sharp(Bw)\,\dd\eta(w)$ for every $B\in\End(V)$;
\item $\lambda_k(A)=m\int w^\sharp((\Lambda^kA)w)\,\dd\eta(w)$ for every $A\in\End(X)$.
\end{enumerate}
\end{theorem}

\begin{proof}
The equivalence of (i) and (ii) is \cref{lem:trace-pairing}, and (i) implies (iii) by taking $B=\Lambda^kA$.

Assume (iii) and put $C=T_\eta-\Id_V$. Then
\[
\tr((\Lambda^kA)C)=0
\qquad(A\in\End(X)).
\]
By \cref{lem:matrix-units}, the operators $\Lambda^kA$ span $\End(V)$, so $\tr(BC)=0$ for every $B\in\End(V)$. Non-degeneracy of the trace pairing gives $C=0$.
\end{proof}

\begin{corollary}\label{cor:auerbach-higher}
Let $(w_i,w_i^\sharp)_{i=1}^m$ be any Auerbach basis of $\Lambda^kX$. Then, for every $A\in\End(X)$,
\begin{equation}\label{eq:auerbach-higher}
\lambda_k(A)=\sum_{i=1}^m w_i^\sharp\big((\Lambda^kA)w_i\big).
\end{equation}
Among isotropic atomic probability measures, the corresponding uniform measure has the smallest possible support.
\end{corollary}

\subsection{The Euclidean Grassmannian formula}

Assume in this subsection that $X=\R^N$ is Euclidean. Let
\[
w_0=e_1\wedge\cdots\wedge e_k.
\]
The orbit of $w_0$ under $O(N)$ is the space of oriented unit simple $k$-vectors and is naturally identified with the oriented Grassmannian
\[
\widetilde G_{k,N}\cong O(N)/(SO(k)\times O(N-k)).
\]
Quotienting by $w\sim-w$ gives the usual Grassmannian
\[
G_{k,N}\cong O(N)/(O(k)\times O(N-k)).
\]
If $E\in G_{k,N}$ and $v_1,\ldots,v_k$ is an orthonormal basis of $E$, then for either orientation
\[
w_E=v_1\wedge\cdots\wedge v_k
\]
one has
\begin{align}\label{eq:grassmann-coefficient}
\ip{(\Lambda^kA)w_E}{w_E}
&=\det\big[\ip{Av_i}{v_j}\big]_{i,j=1}^k
=\det(P_EA\restr_E).
\end{align}
Averaging the rank-one operators $w_E\otimes w_E$ over $G_{k,N}$ gives $m^{-1}\Id_{\Lambda^k\R^N}$, and \eqref{eq:eberlein-intro} follows from \eqref{eq:grassmann-coefficient}. This is Eberlein's formula \cite{Eberlein}; Morrison gave an equivalent formulation on the Stiefel manifold \cite{MorrisonHigher}.

For a general norm, the formulae below integrate over the entire unit sphere of the projective exterior norm. For $2\leqslant k\leqslant N-2$, this sphere contains non-decomposable vectors, so the resulting averages are genuinely different from Grassmannian averages.

\section{Cone measure and polyhedral formulae}\label{sec:cone}

Let $K\subset\R^m$ be a convex body with $0\in\operatorname{int}K$. Its cone probability measure is
\begin{equation}\label{eq:cone-definition}
\nu_K(E)=\frac{\vol\{tx:x\in E,\ 0\leqslant t\leqslant1\}}{\vol(K)}
\qquad(E\subset\partial K\text{ Borel}).
\end{equation}
The ratio is independent of the choice of Lebesgue measure, and for every linear isomorphism $L$ one has $L_\#\nu_K=\nu_{LK}$. At $\mathcal H^{m-1}$-almost every boundary point, the Euclidean outer unit normal $n_K(x)$ exists, and
\begin{equation}\label{eq:cone-density}
\dd\nu_K(x)=\frac{\ip{x}{n_K(x)}}{m\vol(K)}\,\dd\mathcal H^{m-1}(x).
\end{equation}
See \cite{Gardner,Schneider} for standard facts about cone measure and support functions.

\begin{proposition}\label{prop:minkowski-matrix}
For every convex body $K\subset\R^m$,
\begin{equation}\label{eq:minkowski-matrix}
\int_{\partial K}x\otimes n_K(x)\,\dd\mathcal H^{m-1}(x)
=\vol(K)\Id_m.
\end{equation}
Equivalently, for every $B\in\End(\R^m)$,
\begin{equation}\label{eq:minkowski-paired}
\int_{\partial K}\ip{Bx}{n_K(x)}\,\dd\mathcal H^{m-1}(x)
=\tr(B)\vol(K).
\end{equation}
\end{proposition}

\begin{proof}
Convex bodies are sets of finite perimeter, so the Gauss--Green theorem applies; see \cite[Theorem~5.16]{EvansGariepy}. Apply it to the linear vector field $F(x)=Bx$. Since $\div F=\tr(B)$, this gives \eqref{eq:minkowski-paired}. Non-degeneracy of the trace pairing yields \eqref{eq:minkowski-matrix}.
\end{proof}

Now let $V$ be an $m$-dimensional normed space, identify it with $\R^m$ through an auxiliary Euclidean structure, and take $K=B_V$. At a regular point $w\in S_V$, the unique norming functional is
\begin{equation}\label{eq:norming-normal}
w^\sharp(z)=\frac{\ip{z}{n_V(w)}}{\ip{w}{n_V(w)}}.
\end{equation}

\begin{theorem}\label{thm:cone-isotropy}
Cone measure on $S_V$ is isotropic:
\begin{equation}\label{eq:cone-isotropy}
m\int_{S_V}w\otimes w^\sharp\,\dd\nu_V(w)=\Id_V.
\end{equation}
Consequently, for every $B\in\End(V)$,
\begin{equation}\label{eq:cone-trace-general}
\tr(B)
=m\int_{S_V}w^\sharp(Bw)\,\dd\nu_V(w)
=\frac1{\vol(B_V)}\int_{\partial B_V}\ip{Bw}{n_V(w)}\,\dd\mathcal H^{m-1}(w).
\end{equation}
\end{theorem}

\begin{proof}
By \eqref{eq:cone-density} and \eqref{eq:norming-normal},
\begin{align*}
m\int_{S_V}w\otimes w^\sharp\,\dd\nu_V(w)
&=\frac1{\vol(B_V)}
  \int_{\partial B_V}w\otimes n_V(w)\,\dd\mathcal H^{m-1}(w)\\
&=\Id_V
\end{align*}
by \cref{prop:minkowski-matrix}. The trace formula follows from \cref{lem:trace-pairing}.
\end{proof}

\begin{corollary}\label{cor:cone-higher}
Let $V=\Lambda^kX$ with the projective exterior norm and $m=\binom Nk$. Then, for every $A\in\End(X)$,
\begin{align}\label{eq:cone-higher}
\lambda_k(A)
&=m\int_{S_V}w^\sharp\big((\Lambda^kA)w\big)\,\dd\nu_V(w)\notag\\
&=\frac1{\vol(B_V)}
  \int_{\partial B_V}\ip{(\Lambda^kA)w}{n_V(w)}\,\dd\mathcal H^{m-1}(w).
\end{align}
\end{corollary}

The boundary integral becomes a finite sum for polyhedral norms.

\begin{corollary}\label{cor:polyhedral}
Let $V$ be $m$-dimensional and suppose that $B_V$ is a polytope with facets $F_1,\ldots,F_s$, Euclidean outer unit normals $n_1,\ldots,n_s$, and facet centroids
\[
c_j=\frac1{\mathcal H^{m-1}(F_j)}\int_{F_j}w\,\dd\mathcal H^{m-1}(w).
\]
Then, for every $B\in\End(V)$,
\begin{equation}\label{eq:polyhedral-trace}
\tr(B)=\frac1{\vol(B_V)}
\sum_{j=1}^s\mathcal H^{m-1}(F_j)\ip{Bc_j}{n_j}.
\end{equation}
If $V=\Lambda^kX$, then $B$ may be replaced by $\Lambda^kA$ and the left-hand side by $\lambda_k(A)$.
\end{corollary}

\begin{proof}
The normal is constant on the relative interior of each facet, while the lower-dimensional faces have zero $\mathcal H^{m-1}$-measure. Split the final integral in \eqref{eq:cone-trace-general} over the facets and use linearity in $w$.
\end{proof}

\begin{proposition}\label{prop:exterior-polytope}
If $B_X=\operatorname{conv}\{x_1,\ldots,x_r\}$ is a polytope, then the unit ball of the projective exterior norm on $\Lambda^kX$ is the polytope
\begin{equation}\label{eq:exterior-polytope}
B_{\Lambda^kX}
=\operatorname{conv}\{x_{i_1}\wedge\cdots\wedge x_{i_k}:
1\leqslant i_1,\ldots,i_k\leqslant r\}.
\end{equation}
Hence \eqref{eq:polyhedral-trace} gives a finite higher-trace formula whenever $X$ is polyhedral.
\end{proposition}

\begin{proof}
The unit ball of the projective exterior norm is the convex hull of
\[
\{x_1\wedge\cdots\wedge x_k:x_j\in B_X\}.
\]
Writing each $x_j$ as a convex combination of the vertices and expanding by multilinearity shows that every such simple wedge belongs to the right-hand side of \eqref{eq:exterior-polytope}. The reverse inclusion follows directly from the definition of the norm.
\end{proof}

\section{Hypersurface measure and symmetry}\label{sec:symmetry}

Fix an auxiliary Euclidean inner product on an $m$-dimensional normed space $V$. Let
\begin{equation}\label{eq:surface-measure}
\mu_V:=\frac{\mathcal H^{m-1}\restr_{S_V}}{\mathcal H^{m-1}(S_V)}
\end{equation}
be normalised Euclidean hypersurface measure. The norming functional is unique at $\mu_V$-almost every point, by the almost-everywhere differentiability of convex boundaries; see, for example, \cite[Chapter~5]{EvansGariepy}.

\begin{theorem}\label{thm:symmetry}
Let $G$ be a compact subgroup of linear norm isometries of $V$. Assume that the auxiliary Euclidean inner product is $G$-invariant and that
\[
\End(V)^G=\R\Id_V.
\]
Then
\[
T_{\mu_V}=\Id_V.
\]
\end{theorem}

\begin{proof}
Write $T=T_{\mu_V}$. Every $Q\in G$ preserves $S_V$ and $\mu_V$. At every regular point $w$, the norming functional at $Qw$ is
\[
w^\sharp\circ Q^{-1}.
\]
Consequently,
\begin{align*}
QTQ^{-1}
&=m\int_{S_V}(Qw)\otimes(w^\sharp\circ Q^{-1})\,\dd\mu_V(w)\\
&=T.
\end{align*}
Thus $T\in\End(V)^G$, so $T=c\Id_V$. By \eqref{eq:trace-T}, $cm=\tr(T)=m$, and hence $c=1$.
\end{proof}

Every compact group of norm isometries admits an invariant Euclidean inner product, obtained by averaging any inner product over Haar probability measure. The corresponding hypersurface measure is therefore naturally adapted to the chosen symmetry group.

For a finite orthogonal group $G$, the scalar-commutant condition is equivalent to the twirling identity
\begin{equation}\label{eq:twirl}
\frac1{|G|}\sum_{Q\in G}Q^{-1}BQ
=\frac{\tr(B)}m\Id_V
\qquad(B\in\End(V)).
\end{equation}
Indeed, the left-hand side is the Hilbert--Schmidt orthogonal projection of $\End(V)$ onto $\End(V)^G$. In particular, for every Euclidean vector $v$,
\begin{equation}\label{eq:orbit-second-moment}
\frac1{|G|}\sum_{Q\in G}(Qv)\otimes(Qv)
=\frac{\norm{v}_2^2}{m}\Id_V.
\end{equation}
Thus every orbit satisfies the homogeneous second-moment identity characteristic of a spherical $2$-design.

We now apply this to exterior powers. Let $e_1,\ldots,e_N$ be a basis of $X$, declare it Euclidean orthonormal, and let $\mathcal B_N$ be the group of signed permutation matrices.

\begin{lemma}\label{lem:hyperoctahedral}
For every $1\leqslant k\leqslant N$,
\[
\End(\Lambda^k\R^N)^{\Lambda^k\mathcal B_N}
=\R\Id_{\Lambda^k\R^N}.
\]
\end{lemma}

\begin{proof}
For a $k$-element set $I\subset\{1,\ldots,N\}$, let $e_I$ denote the corresponding exterior basis vector. A diagonal sign change $D_\varepsilon$ acts by
\[
(\Lambda^kD_\varepsilon)e_I
=\left(\prod_{i\in I}\varepsilon_i\right)e_I.
\]
Distinct $k$-element subsets give distinct characters of the sign-change subgroup. Hence any operator commuting with all sign changes is diagonal in the basis $(e_I)$. Permutation matrices act transitively on the $k$-element subsets, so commuting with all permutations forces all diagonal entries to coincide.
\end{proof}

\begin{corollary}\label{cor:one-symmetric}
Suppose that $X$ has a $1$-symmetric basis. For every $1\leqslant k\leqslant N$, normalised Euclidean hypersurface measure on the unit sphere of $\Lambda^kX$, formed from the Euclidean structure for which that basis is orthonormal, is isotropic. Therefore
\begin{equation}\label{eq:surface-higher}
\lambda_k(A)
=\binom Nk\int_{S_{\Lambda^kX}}
 w^\sharp\big((\Lambda^kA)w\big)\,\dd\mu_{\Lambda^kX}(w)
\end{equation}
for every $A\in\End(X)$.
\end{corollary}

\begin{proof}
The signed permutation group acts by norm isometries on $X$, hence by \cref{lem:exterior-bound} on every $\Lambda^kX$. Apply \cref{thm:symmetry,lem:hyperoctahedral}.
\end{proof}

The same argument applies to many smaller groups. For example, the dihedral group $D_q$ with $q\geqslant3$ suffices in dimension two, and the symmetry groups of the Platonic solids suffice on $\R^3$ and, through the Hodge isomorphism, on $\Lambda^2\R^3$. The scalar-commutant condition, rather than the size of the group, is the essential point.

\section{Power-weighted measures and first variations}\label{sec:variation}

Let $m\geqslant2$, and let $h\in C^2(\sphere^{m-1})$ be the support function of a strictly convex body $K_h$ with $C^2$ boundary and positive Gauss curvature. Thus
\[
L_h(u):=\nabla_S^2h(u)+h(u)\Id_{T_u\sphere^{m-1}}
\]
is positive definite. The inverse Gauss map is
\begin{equation}\label{eq:gauss-param}
x_h(u)=h(u)u+\nabla_Sh(u),
\end{equation}
and the surface-area element is
\begin{equation}\label{eq:surface-density-support}
\dd\mathcal H^{m-1}(x_h(u))
=\abs{\sphere^{m-1}}\det L_h(u)\,\dd\omega(u),
\end{equation}
where $\omega$ is rotationally invariant probability measure. The norming functional at $x_h(u)$ is represented by $u/h(u)$.

For $\alpha\in\R$, define
\begin{equation}\label{eq:Z-alpha}
Z_\alpha(h):=\int_{\sphere^{m-1}}h(u)^\alpha\det L_h(u)\,\dd\omega(u)
\end{equation}
and let $\nu_h^\alpha$ be the probability measure whose pullback under \eqref{eq:gauss-param} is
\begin{equation}\label{eq:alpha-measure}
\frac{h(u)^\alpha\det L_h(u)}{Z_\alpha(h)}\,\dd\omega(u).
\end{equation}
Thus $\nu_h^0$ is normalised hypersurface measure and $\nu_h^1$ is cone measure. Its isotropy tensor is
\begin{equation}\label{eq:T-alpha}
T_\alpha(h)
=\frac{m}{Z_\alpha(h)}
\int_{\sphere^{m-1}}x_h(u)\otimes u\,
 h(u)^{\alpha-1}\det L_h(u)\,\dd\omega(u),
\end{equation}
where $u$ is identified with the Euclidean functional $z\mapsto\ip{z}{u}$.

\begin{theorem}\label{thm:first-variation}
Let $g\in C^\infty(\sphere^{m-1})$ and $h_\varepsilon=1+\varepsilon g$. For all sufficiently small $\varepsilon$, assume that $h_\varepsilon$ is a support function as above. Then
\begin{equation}\label{eq:first-variation-operator}
\left.\frac{\dd}{\dd\varepsilon}\right|_{\varepsilon=0}T_\alpha(h_\varepsilon)
=m(\alpha-1)\int_{\sphere^{m-1}}g(u)
\left(u\otimes u-\frac1m\Id_m\right)\,\dd\omega(u).
\end{equation}
Equivalently, for every $B\in\End(\R^m)$,
\begin{equation}\label{eq:first-variation-paired}
\left.\frac{\dd}{\dd\varepsilon}\right|_{\varepsilon=0}
\tr\big(BT_\alpha(h_\varepsilon)\big)
=m(\alpha-1)\int_{\sphere^{m-1}}g(u)
\left(\ip{Bu}{u}-\frac{\tr(B)}m\right)\,\dd\omega(u).
\end{equation}
\end{theorem}

\begin{proof}
Set
\[
f_B(u)=\ip{Bu}{u},
\qquad
\overline g=\int_{\sphere^{m-1}}g\,\dd\omega.
\]
The first-order expansions are
\begin{align*}
x_{h_\varepsilon}(u)
&=u+\varepsilon\big(g(u)u+\nabla_Sg(u)\big),\\
h_\varepsilon^{\alpha-1}
&=1+\varepsilon(\alpha-1)g+O(\varepsilon^2),\\
\det L_{h_\varepsilon}
&=1+\varepsilon\big(\Delta_Sg+(m-1)g\big)+O(\varepsilon^2).
\end{align*}
Consequently,
\begin{equation}\label{eq:Z-derivative}
Z_\alpha(h_\varepsilon)
=1+\varepsilon(\alpha+m-1)\overline g+O(\varepsilon^2).
\end{equation}
Pairing the numerator in \eqref{eq:T-alpha} with $B$, differentiating at zero, and dividing by $m$ gives
\begin{equation}\label{eq:numerator-derivative}
\int_{\sphere^{m-1}}
\left[
\ip{B\nabla_Sg}{u}
+(\alpha+m-1)gf_B
+(\Delta_Sg)f_B
\right]\,\dd\omega.
\end{equation}
We use two elementary identities on the sphere. First,
\[
\div_S\big(P_{T_u\sphere^{m-1}}B^\top u\big)
=\tr(B)-mf_B(u),
\]
so integration by parts gives
\begin{equation}\label{eq:sphere-ibp-1}
\int\ip{B\nabla_Sg}{u}\,\dd\omega
=m\int gf_B\,\dd\omega-\tr(B)\overline g.
\end{equation}
Second,
\[
\Delta_Sf_B=-2m\left(f_B-\frac{\tr(B)}m\right),
\]
and hence
\begin{equation}\label{eq:sphere-ibp-2}
\int(\Delta_Sg)f_B\,\dd\omega
=-2m\int gf_B\,\dd\omega+2\tr(B)\overline g.
\end{equation}
Combining \eqref{eq:numerator-derivative}--\eqref{eq:sphere-ibp-2}, the derivative of the numerator in \eqref{eq:T-alpha}, after multiplication by $m$, is
\[
m(\alpha-1)\int gf_B\,\dd\omega+m\tr(B)\overline g.
\]
The derivative of the quotient also contains the normalising contribution from \eqref{eq:Z-derivative}. Since the value at $\varepsilon=0$ is $\tr(B)$, the resulting derivative is
\begin{align*}
&m(\alpha-1)\int gf_B\,\dd\omega
+m\tr(B)\overline g
-(\alpha+m-1)\tr(B)\overline g\\
&\hspace{4em}=m(\alpha-1)\int g\left(f_B-\frac{\tr(B)}m\right)\,\dd\omega.
\end{align*}
This proves \eqref{eq:first-variation-paired}; non-degeneracy of the trace pairing gives \eqref{eq:first-variation-operator}.
\end{proof}

Each matrix entry of
\[
u\longmapsto u\otimes u-\frac1m\Id_m
\]
is a degree-$2$ spherical harmonic. Therefore the integral in \eqref{eq:first-variation-operator} depends only on $P_2g$, the orthogonal projection of $g$ onto the degree-$2$ harmonic subspace. More explicitly, the map
\[
\Sym_0(\R^m)\longrightarrow\mathcal H_2,
\qquad
S\longmapsto\big[u\mapsto u^\top Su\big],
\]
is an $O(m)$-equivariant isomorphism. The fourth-moment identity
\begin{equation}\label{eq:fourth-moment}
\int_{\sphere^{m-1}}u_iu_ju_ku_\ell\,\dd\omega(u)
=\frac{\delta_{ij}\delta_{k\ell}+\delta_{ik}\delta_{j\ell}+\delta_{i\ell}\delta_{jk}}{m(m+2)}
\end{equation}
then gives the following consequence; see \cite{DaiXu} for background on spherical harmonics.

\begin{corollary}\label{cor:degree-two}
If
\[
P_2g(u)=u^\top Su,
\qquad S\in\Sym_0(\R^m),
\]
then
\begin{equation}\label{eq:quadratic-variation}
\left.\frac{\dd}{\dd\varepsilon}\right|_{\varepsilon=0}T_\alpha(h_\varepsilon)
=\frac{2(\alpha-1)}{m+2}S.
\end{equation}
In particular, $\alpha=1$ is the unique real exponent such that $T_\alpha(h)=\Id_m$ for every smooth strictly convex centrally symmetric body containing the origin.
\end{corollary}

\begin{proof}
By \eqref{eq:fourth-moment},
\[
\int_{\sphere^{m-1}}(u^\top Su)
\left(u\otimes u-\frac1m\Id_m\right)\,\dd\omega(u)
=\frac{2}{m(m+2)}S.
\]
Substitute this into \eqref{eq:first-variation-operator}. Cone measure, corresponding to $\alpha=1$, is isotropic for every convex body by \cref{thm:cone-isotropy}. If $\alpha\neq1$, choose $S\neq0$ and $g(u)=u^\top Su$; then \eqref{eq:quadratic-variation} is non-zero, so $T_\alpha(h_\varepsilon)\neq\Id_m$ for all sufficiently small non-zero $\varepsilon$.
\end{proof}

For hypersurface measure, $\alpha=0$, and \eqref{eq:first-variation-operator} becomes
\begin{equation}\label{eq:hypersurface-variation}
T_0(h_\varepsilon)
=\Id_m-m\varepsilon\int_{\sphere^{m-1}}g(u)
\left(u\otimes u-\frac1m\Id_m\right)\,\dd\omega(u)
+O(\varepsilon^2).
\end{equation}
Thus the degree-$2$ component is the complete first-order obstruction.

\subsection{Explicit failures of hypersurface isotropy}

\begin{example}\label{ex:hexagon}
On $\R^2$, let
\[
\norm{(x,y)}=\max\{|x|,|y|,|y-x|\}.
\]
The unit ball is the hexagon with vertices
\[
(1,1),(0,1),(-1,0),(-1,-1),(0,-1),(1,0).
\]
Using normalised arclength measure and the constant norming functional on the relative interior of each edge gives
\begin{equation}\label{eq:hexagon-tensor}
T_{\mu}
=\begin{pmatrix}
1&2-\dfrac{3\sqrt2}{2}\\[3pt]
2-\dfrac{3\sqrt2}{2}&1
\end{pmatrix}\neq\Id_2.
\end{equation}
For
\[
A=\begin{pmatrix}0&1\\0&0\end{pmatrix},
\]
one has $\tr(A)=0$ but
\[
\tr(AT_\mu)=2-\frac{3\sqrt2}{2}\neq0.
\]
\end{example}

\begin{example}\label{ex:smooth}
In $\R^2$, take
\[
h_\varepsilon(\theta)=1+\varepsilon\cos(2\theta).
\]
Formula \eqref{eq:hypersurface-variation} gives
\begin{equation}\label{eq:smooth-example}
T_0(h_\varepsilon)
=\Id_2+\frac\varepsilon2
\begin{pmatrix}-1&0\\0&1\end{pmatrix}
+O(\varepsilon^2).
\end{equation}
For $A=\operatorname{diag}(1,-1)$,
\[
\tr(AT_0(h_\varepsilon))-\tr(A)
=-\varepsilon+O(\varepsilon^2).
\]
\end{example}

Although the first variation in \eqref{eq:hypersurface-variation} is symmetric, the full hypersurface tensor need not be. In dimension two, write
\[
u(\theta)=(\cos\theta,\sin\theta),
\qquad
t(\theta)=(-\sin\theta,\cos\theta).
\]
For a $C^2$ support function $h$ with $h+h''>0$, a direct use of
\[
x=h u+h't,
\qquad
\dd s=(h+h'')\dd\theta,
\qquad
x^\sharp=\frac{u}{h}
\]
gives
\begin{equation}\label{eq:skew-exact}
T_0(h)-T_0(h)^\top
=\frac{2}{\int_0^{2\pi}h\,\dd\theta}
\left(\int_0^{2\pi}\frac{h'h''}{h}\,\dd\theta\right)
\begin{pmatrix}0&-1\\1&0\end{pmatrix}.
\end{equation}
For example, if
\[
h(\theta)=1+\varepsilon\big(\cos2\theta+\sin4\theta\big),
\]
then
\begin{equation}\label{eq:skew-example}
T_0(h)-T_0(h)^\top
=12\varepsilon^3
\begin{pmatrix}0&1\\-1&0\end{pmatrix}
+O(\varepsilon^4).
\end{equation}
Thus symmetry of the first-order term is a linearised phenomenon rather than a general property of hypersurface averages.

\section{Concluding remarks}

The isotropy criterion separates three rather different constructions. Auerbach bases provide optimal finite support but are highly non-canonical. Cone measure is canonical and works for every norm. Hypersurface measure is geometrically familiar but generally fails without an appropriate symmetry group.

For $k\geqslant2$, the cone and hypersurface formulae on $\Lambda^kX$ integrate over the whole unit sphere of the projective exterior norm. A separate problem is to determine when an isotropic measure can be chosen on the subset of decomposable unit vectors, thereby retaining a direct interpretation in terms of $k$-dimensional subspaces. In the Euclidean case this is exactly the Grassmannian formula \eqref{eq:eberlein-intro}.

Another natural problem is to characterise those norm balls whose normalised Euclidean hypersurface measure is isotropic. \Cref{thm:symmetry} supplies a broad symmetry-based sufficient condition, while \cref{thm:first-variation,cor:degree-two} identify the complete infinitesimal obstruction at the Euclidean ball.

\paragraph*{Acknowledgements.}
I thank Kent Morrison for drawing my attention to Eberlein's higher-trace formulae and for sharing his unpublished note \cite{MorrisonHigher}. I also thank Anthony Quas for helpful discussions about exterior powers of operators on finite-dimensional normed spaces.

\end{document}